\documentclass[amsa]{ipart}

\newtheorem{alg}{\bf Algorithm}
\newtheorem{theorem}{\bf Theorem}
\newtheorem{exm}{\bf Example}

\newtheorem{lem}{\bf Lemma}

\usepackage{graphicx}
\usepackage{amssymb}

\begin{document}

\begin{frontmatter}


\title{An Inexact Inverse Power Method for Numerical Analysis of Stochastic Dynamic Systems}
\begin{aug}
\author{\fnms{Yuquan Sun} \thanksref{t1}\ead[label=e1]{sunyq@buaa.edu.cn}},
\address{LMIB \& School of Mathematics and Systems Science,\\
 BeiHang University, Beijing China, 100191.\\
\printead{e1}}
\author{\fnms{Fanghui Gong}, \ead[label=e2]{fanghui.gong@buaa.edu.cn}}
\address{LMIB \& School of Mathematics and Systems Science,\\
 BeiHang University, Beijing China, 100191.\\
\printead{e2}}

\author{\fnms{Igor V. Ovchinnikov}, \ead[label=e3]{iovchinnikov@ucla.edu}}
\address{Electrical Engineering Department, \\
University of California at Los Angeles, Los Angeles, CA 90095.\\
\printead{e3}}
\author{\fnms{Kang L. Wang} \thanksref{t2}\ead[label=e4]{wang@ee.ucla.edu}}
\address{Electrical Engineering Department, \\
University of California at Los Angeles, Los Angeles, CA 90095.\\
\printead{e4}}

\thankstext{t1}{Thanks to National Science Foundation of China (No. 11201020 )for support.}
\thankstext{t2}{Thanks to Raytheon Endowed Professorship for support.}
\runauthor{Yuquan Sun et al.}

\end{aug}



\begin{abstract}
This paper proposes an efficient method for computing partial eigenvalues of large sparse matrices what can be called the inexact inverse power method (IIPM). It is similar to the inexact Rayleigh quotient method and inexact Jacobi-Davidson method that it uses only a low precision approximate solution for the inner iteration. But this method uses less memory than inexact Jacobi-Davidson method and has stronger convergence performance than inexact Rayleigh quotient method. We exemplify the advantages of IIPM by applying it to find the ground state in theory of stochastics. Here we need to solve hundreds of large-scale matrix. The computational results show that this approach is a particularly useful method.
\end{abstract}

\end{frontmatter}


\section{Introduction} 
Stochastic differential equations (SDEs, see, \emph{e.g.}, Refs. \cite{DOTO} and Refs. therein) is a class of mathematical models with the widest applicability in modern science.
\begin{eqnarray}
\dot x(t) = F(x(t)) + (2\Theta)^{1/2}e_{a}(x(t))\xi^a(t),
\end{eqnarray}
where $x\in X$ is a point in the phase space $X$, which is a topological manifold, $F\in TX$ is a vector field from the tangent space of $X$, $TX$, called the flow vector field, $\Theta$ is the intensity of temperature of the Gaussian white noise, $\xi$, with the standard expectation values $\langle \xi^a(t) \rangle =0,  \langle \xi^a(t)\xi^b(t') \rangle = \delta(t-t')\delta^{ab}$, and the $a$-set of vector fields $e_a\in TX$ defining the coupling of the noise to the system.

In physics, for example, it describes everything in nature above the scale of quantum degeneracy/coherence. One of the main statistical characteristics of these systems is the probability density of the solution
of this equation. The probability density can be studied by the corresponding Fokker-Planck equation. It can be perform further studies by the supersymmetric theory of stochastics (STS) \cite{Igor1,Igor1_1,Igor2}, which is one of the latest advancements in the theory of SDEs. Among a few other important findings, STS seems to explain 1/f noise \cite{DOTO}, power-laws statistics of various avalanche-type processes \cite{Aschwanden} and other realizations of mysterious and ubiquitous dynamical long-range order \cite{Igor1} in nature.

As compared to classical approaches to SDEs, STS differs in two fundamental ways. First, the Hilbert space of a stochastic model in STS, $\mathcal{H}$, is the entire exterior algebra of the phase space, \emph{i.e.}, the space of differential forms or $k$-forms of all degrees,
\begin{eqnarray}
\psi^{(k)} =
\psi^{(k)}_{i_1...i_k}(x)dx^{i_1}\wedge...\wedge dx^{i_k}\in \Omega^{k},
\mathcal{H} = \bigoplus\nolimits_{k=0}^{D}\Omega^{k},
\end{eqnarray}
where $\psi_{i_1...i_k}$ is an antisymmetric contravariant tensor, $\Omega^{k}$ is the space of all $k$-forms and $D$ is the dimensionality of the phase space of the model. This picture generalizes the classical approach to SDEs, where the Hilbert space is thought to be the space of only top differential forms that have the meaning of total probability distributions in a coordinate-free setting.

The second distinct feature of the STS, is that the finite-time stochastic evolution operator (SEO) has a clear mathematical meaning. Specifically,
\begin{eqnarray}
\psi(t) = \hat{\mathcal{M}}_{tt'}\psi(t'), \hat{\mathcal{M}}_{tt'} = \langle M^*_{t't}\rangle,
\end{eqnarray}
where $M^*_{t't}$ is a pullback or action induced by the SDE-defined noise-configuration-dependent diffeomorphism $M_{t't}$, so that a noise-configuration-dependent solution of SDE with initial condition $x(t)|_{t=t_0}=x_0$ can be given as $x(t) = M_{tt_0}(x_0)$, and brackets denote stochastic averaging over all the configuration of the noise.

The finite-time SEO can be shown \cite{Igor1} to be,
\begin{eqnarray}
\hat{\mathcal{M}}_{tt'} = e^{-(t-t')\hat H},\label{F_T_SEO}
\end{eqnarray}
where the (infinitesemal) SEO is given as,
\begin{eqnarray}
\hat H = \hat {\mathcal L}_{F} - \Theta \hat {\mathcal L}_{e_a}\hat {\mathcal L}_{e_a},\label{Lieseo}
\end{eqnarray}
where $\hat {\mathcal L}$ is a Lie or physical derivative along the corresponding vector field.

The presence of the topological supersymmetry given by Eq.(\ref{F_T_SEO}) tailors the following properties of the eigensystem of the SEO. Here are two types of eigenstates. The first type is the supersymmetric singlets that are non-trivial in the De Rahm cohomology. Each De Rahm cohomology class of $X$ must provide one supersymmetric singlet \cite{Ours}. All supersymmetic eigenstates have exactly zero eigenvalue.
The second type of state are non-supersymmetric doublets.
There are no restrictions on the eigenvalues of the non-supersymemtric eignestates other than they must be either real or come in complex conjugate pairs known in the dynamical systems theory as Ruelle-Pollicott resonances and that the real part of its eigenvalue must be bounded from below in case when the diffusion part of the SEO is elliptic. Most of the eigenstates of the SEO are non-supersymmetric. In particular, all eigenstates with non-zero eigenvalues are non-supersymmetric.

The ground state(s) is the state(s) with the lowest real part of its eigenvalue. As is seen from the exponential temporal evolution in Eq.(\ref{F_T_SEO}), the ground state(s) grows (and oscillates if its eigenvalue is complex) faster than any other eigenstate. When the ground state is a non-supersymmetric eigenstate, it is said that the topological supersymmetry is spontaneous broken. The topological supersymmetry breakdown can be identified with the stochastic generalization of the concept of deterministic chaos \cite{Igor1,Igor2}, and this identification is an important finding for applications.

Whether the topological supersymmetry is spontaneously broken or not can be unambiguously determined from the eigensystem of the SEO.
So the numerical investigation of the SEO's eigensystem is an important method. Because different parameters will give different eigensystems,
we need to solve hundreds of eigenvalue problems.

Eigenvalue problem of sparse matrices is an important problem that has applications in many branches of modern science. This problem has a long history and several powerful methods of the numerical studies of sparse matrices have been proposed and implemented by now. One of the most successful such implementations is ARPACK \cite{Sorensen2} based on the Implicitly Restarted Arnoldi Method \cite{Sorensen}. ARPACK is a collection of Fortran subroutines designed to compute a few eigenvalues and corresponding eigenvectors of a sparse matrix and it is the foundation of the commonly used MATLAB command "eigs".

In many applications, one has to compute eigenvalues with the smallest real part, \emph{i.e.}, the leftmost in the complex plane. On the other hand, the structure of the Arnoldi Method targets eigenvalues with the largest magnitude. Therefore, for "low-lying" eigenvalues, the "eigs" function may encounter convergence problems, even when using a large trial subspace.

The problem of low-lying eigenvalue is better addressed with the inverse power method that transforms it into the largest eigenvalue problem. Yet another generalization is the Shift-Invert Arnoldi method, which is the original Arnoldi method applied to the shift-inverted matrix $B = (A - \sigma I)^{-1}$, so that it can find eigenvalues near to the given target $\sigma$. There exist other variations of the parental Arnoldi method including the Residual Arnoldi and the Shift-Invert Residual Arnoldi methods \cite{Lee}.

One of the problems or the Shift-Invert Arnoldi method is that the inverse matrix $(A - \sigma I)^{-1}$ cannot be easily computed for large matrices. This inversion is practically achieved by iteratively solving the corresponding system of linear equations (CSLE). This may already be a difficult problem for large matrices.

Yet another approach is the use of inexact methods \cite{Freitag, Simoncini}. The main idea of these methods is
computing an approximate solution of the inner equation. The convergence analysis of the
inexact method has already been widely studied \cite{Notay, Xue}. Recently, a general convergence theory of the Shift-Invert Residual Arnoldi (SIRA) method has been established \cite{Jia}.

In order to ensure good convergence, these methods need to expand the dimensionality of the working subspace continuously from iteration to iteration. For large problems, the computation and storage costs may be very high. As it turns out, in our application, we need to solve hundreds of large matrix, under limited time and resource constraints. In order to achieve this goal, we propose what we call the inexact inverse power method (IIPM). The advantages are that one only needs to store two vectors (two-dimensional subspace) during all iterations and for this reason save considerably the required computation resources. At the same time, it keeps a high convergence rate. The existing convergence analyses are based on the prior knowledge of the eigenvalue information. In theory, the convergence of this method can be guaranteed, but lack practical guidance for real computation.
From a view of ensuring the convergence, we analyze the convergence of the new algorithm and propose a convergence criterion of inner iteration for practical computation.

The paper is organized as follows. In Section \ref{s3}, we describe the proposed IIPM and  analyze the convergence. In Section \ref{s5}, we exemplify the advantages of the IIMP by applying it to the problem of the diagonalization of the stochastic evolution operators of the ABC and Kuramoto models. Section \ref{conclusion} concludes this paper.


\section{The Inexact Inverse Power Method}
\label{s3}

In this section we would like to discuss the theory of the IIPM for the large matrix diagonalization problems. As we mentioned in the Introduction, this method is derivative of its parental IPM. Therefore, we begin the discussion with the introduction of the IPM.

\begin{alg} {\bf Inverse power method}\label{Ipower}
\begin{description}
  \item[1:]Given starting vector $x_1$, and convergence criterion  $tol$  \\[-5mm]
  \item[2:] for $i = 1,2,\ldots,n$                      \\[-5mm]
  \item[3:]\hspace{0.5cm} $y = A^{-1}x_i$               \\[-5mm]
  \item[4:]\hspace{0.5cm} $x_{i+1} = y/\|y\|$           \\[-5mm]
  \item[5:]\hspace{0.5cm} $\lambda = x_{i+1}^{T} A x_{i+1}$  \\[-5mm]
  \item[6:]\hspace{0.5cm} $t = Ax_{i+1}- \lambda x_{i+1}$  \\[-5mm]
  \item[7:]\hspace{0.5cm}   if $ \|t\|\leq tol $, break \\[-5mm]
  \item[8:] end for                                     \\[-5mm]
\end{description}
\end{alg}

We can apply this process to matrix $(A-\sigma I)^{-1}$ instead of $A$, where $\sigma$ is called a shift. This will allow us to compute the eigenvalue closest to $\sigma$. When $\sigma$ is very close to the desired eigenvalue, we can obtain a faster convergence rate.

For large scale matrix, neither $y = A^{-1}x_i$ nor $y =(A-\sigma I)^{-1}x_i$ can be computed directly. It is difficult to obtain an accurate solution even by solving the corresponding linear systems,
\begin{equation}
(A-\sigma I)y =x_i.\label{equ2}
\end{equation}

According to the idea of inexact methods, we can use an iterative approach to compute an approximate solution $\tilde{y}$ of Eq.(\ref{equ2}).
Namely, we can use ${x}_{i+1} = \tilde{y}/\|\tilde{y}\|$ as the updated approximate eigenvector. This alteration of the IPM leads one to the IIPM.

\begin{alg} {\bf  Inexact inverse  power method}\label{vpower}
\begin{description}
  \item[1:]Given a target $\sigma $, starting vector $x_1$ and  convergence criterion  $\varepsilon_1,\varepsilon_2$,   \\[-5mm]
  \item[2:] for $k = 1,2,\ldots,n$                                                                        \\[-5mm]
  \item[3:]\hspace{0.5cm} Compute an approximate solution $\tilde{y}$ of $(A-\sigma I )y = x_i$  with
                $$\| r \|=\| x_i - (A-\sigma I )\tilde{y} \| < \varepsilon_1 $$
 \item[4:]\hspace{0.5cm}  Compute  eigenpairs$(\lambda,x_{i+1}) $ from span\{$x_i, \tilde{y}\}$     \\[-5mm]
  \item[5:]\hspace{0.5cm} $t = Ax_{i+1}- \lambda x_{i+1}$                                                \\[-5mm]
  \item[6:]\hspace{0.5cm} if $\|t\| \leq \varepsilon_2$, break                                             \\[-5mm]
  \item[7:]\hspace{0.5cm} end for                                                                          \\[-5mm]
 \end{description}
\end{alg}

When we use the inexact solution $\tilde{y}$ instead of the exact solution $y$, the convergence property of
${{x}}_{i+1}$ is the most important issue. This means that we need a quantitative standard for $\varepsilon $.
In the kind of inexact methods, the convergence is obtained  by
 analyzing  the ability of $\tilde{y}$ to mimic $y$ and the convergence is guaranteed by the subspace expanding.
For this method, we write the approximate solution ${\tilde{y}}$ of  Eq.(\ref{equ2}) as the exact solution of the following  perturbed equation
\begin{equation}
(A-\sigma I + \delta A )\tilde{y} =x_i,\label{equ3}
\end{equation}
here $\delta A$ is the perturbation matrix of $(A-\sigma I)$. The residual of  Eq.(\ref{equ2}) can be written as
$ r = \delta A \tilde{y} $.

\begin{lem} \label{lem1}  The approximate solution ${\tilde{y}}$ and the exact solution $y$ of  Eq.(\ref{equ2}) have the following relationship
\begin{equation}
 y - \tilde{y} \approx (A-\sigma I)^{-1} \delta A y.
\end{equation}

{\bf Proof:} For a matrix $X$ and the corresponding unit matrix $I$, if $\|X\| < 1$, then $I-X$ is invertible \cite{Demmel} and
\begin{equation}
(I-X)^{-1} = \sum_{i=0}^\infty X^i .\label{inverse}
\end{equation}

Now, we can get
\begin{equation}
\begin{array}{ccl}
\tilde{y}& = & (A-\sigma I + \delta A )^{-1} x_i \\
         & = & (I + (A-\sigma I)^{-1} \delta A )^{-1} (A-\sigma I)^{-1} x_i.
\end{array}
\end{equation}
We use formula (\ref{inverse}) to $(I + (A-\sigma I)^{-1} \delta A )^{-1}$ and ignore higher order terms to obtain
\begin{equation}
\begin{array}{ccl}
\tilde{y}& \approx & [I - (A-\sigma I)^{-1} \delta A ] (A-\sigma I)^{-1} x_i \\
         & = & [I - (A-\sigma I)^{-1} \delta A ]y                             \\
         & = & y - (A-\sigma I)^{-1} \delta A y.
\end{array}
\end{equation}
\end{lem}
Then we obtain the result about $\tilde{y}$ and $y$.

Suppose $(\lambda ,x) $ is a simple desired eigenpair of $A$ and $(\frac{1}{\lambda - \sigma},x)$
is a simple eigenpair of $(A-\sigma I)^{-1}$. In Algorithm \ref{Ipower}, both  $y$ and ${x}_{i}$ are approximate eigenvectors but $y$ is a better approximate eigenvector than ${x}_{i}$. This can be obtained from the convergence properties of the power method.
In Algrithm 2, it is hard to ensure that $\tilde{y}$ is better than ${x}_{i}$  only by the relationship between $\tilde{y}$ and $y$.
To obtain the convergence property of Method \ref{vpower}, we need an $ \varepsilon $  to ensure
\begin{equation}
\|{x}_{i+1}-x\|\leq \|{x}_{i}-x\|,  \label{relation}
\end{equation}
here ${x}_{i+1} = \tilde{y}/ \|\tilde{y}\|.$

Since $x,{x}_{i}, {x}_{i+1}$ are unit vectors, their relationship can be expressed better by the angle.
The relationship (\ref{relation}) is equivalent to
\begin{equation}
\mbox{sin} \angle ({x}_{i+1},x) \leq \mbox{sin} \angle({x}_{i},x).  \label{relation2}
\end{equation}

For the convenience of analysis, we set
$$B = (A-\sigma I)^{-1}.$$
The eigenvalues of $B$ satisfy $\mu_1 \gg \mu_2\geq \mu_3\geq\cdots \geq \mu_n$. $ (\mu_1,x)$ is the desired eigenpair of $B$.
Let $(x,X_{\bot})$ be a unitary matrix, where $\mbox{span}\{ X_{\bot}\}$ is the orthogonal complement of ${x}$.
Then $x_i,y$ and $\delta A y$ can be expressed as
\begin{equation}
\begin{array}{ccl}
x_i       & = & \alpha x + \beta z,\\
y         & = & \mu_1 \alpha x +  \beta Bz,\\
\delta A y& = & \tilde{\alpha} x + \tilde{\beta} \tilde{z},
\end{array}\label{relation1}
\end{equation}
where $z,\tilde{z} \in \mbox{span}\{ X_{\bot}\}$, $\|\delta A \|\leq \varepsilon$.

By Lemma \ref{lem1}, $\tilde{y}$ can be written as
\begin{equation}
\tilde{y} =\mu_1 \alpha x +  \beta B z +
\mu_1 \tilde{\alpha} x +  \tilde{ \beta} B \tilde{z}.
\end{equation}

For the convergence of Algorithm \ref{vpower}, we have the following result.
\begin{theorem} Suppose $B$ is symmetric, $x$ is the desired eigenvector, $x_i$ is the current approximation of $x$.
 $y$ and $\tilde{y}$ are exact and approximate solution of  Eq.(\ref{equ2}).
If $\varepsilon $ satisfies
 $$\varepsilon < \frac{(1- 2\mu_2/\mu_1)\beta\alpha}{(2\mu_2/\mu_1) \alpha + \beta}.$$
Then we have
$$\mbox{sin}\angle(\tilde{y},x) < \mbox{sin}\angle(x_i,x).$$
{\bf Proof:}
Eq. (\ref{relation1}) shows that $\mbox{tan}\angle(x_i,x)= \frac{|\beta|}{|\alpha|}$.
$\mbox{tan}\angle(y,x)= \frac{|\beta|\|Bz\|}{|\alpha \mu_1|}$.
Since $z \in \mbox{span}\{ X_{\bot}\}$, so $\|Bz\| \leq \mu_2$. Then we have
$$\mbox{tan}\angle(y,x)< \frac{|\beta \mu_2|}{|\alpha \mu_1|} < \frac{\mu_2}{\mu_1}\mbox{tan}\angle(x_i,x).$$
We can write $\tilde{y}$ as $\tilde{y} =\mu_1( \alpha + \tilde{\alpha}) x +  (\beta + \tilde{\beta}) B(z + \tilde{z} )$.
From $z,\tilde{z} \in \mbox{span}\{ X_{\bot}\}$, we obtain $\|B(z + \tilde{z}) \|\leq 2 \mu_2$.
From $|\tilde{\beta}|< \varepsilon, |\tilde{\alpha}|< \varepsilon$, we get $|\beta + \tilde{\beta}| \leq |\beta| + \varepsilon$
and $|\alpha + \tilde{\alpha}| \geq  |\alpha| - \varepsilon$.
For the angle between $\tilde{y}$ and $x$,  we have inequality

\begin{equation}\label{speed}
\mbox{tan}\angle(\tilde{y},x)\leq \frac{2 \mu_2 |\beta + \tilde{\beta}|}{\mu_1|\alpha + \tilde{\alpha}|}
\leq \frac{2 \mu_2}{\mu_1} \frac{|\beta| + \varepsilon}{|\alpha| - \varepsilon}.
\end{equation}

If
\begin{equation}
\frac{2 \mu_2}{\mu_1} \frac{|\beta| + \varepsilon}{|\alpha| - \varepsilon} < \frac{|\beta|}{|\alpha|},
\label{iequal2}
\end{equation}
then we obtain $\mbox{tan}\angle(\tilde{y},x) < \mbox{tan}\angle(x_i,x)$.

Since (\ref{iequal2}) is equivalent to
\begin{equation}
 \varepsilon < \frac{(1- (2\mu_2/\mu_1))|\beta\alpha|}{(2\mu_2/\mu_1) |\alpha| + |\beta|},
\label{iequal3}
\end{equation}
then we finish the proof.
\end{theorem}

When the angle between $x_i$ and $x$ is not very small, the values of $|\alpha|$ and $|\beta|$ have same order $|\alpha|=O(|\beta|)$.
From $\mu_1 \gg \mu_2$, we obtain  $(2\mu_2/\mu_1)$ is a small number and $1- (2\mu_2/\mu_1)\approx 1$. So the requirements of
$ \varepsilon$ is
$$ \varepsilon < \frac{(1- (2\mu_2/\mu_1))|\beta\alpha|}{(2\mu_2/\mu_1) |\alpha| + |\beta|}= O( |\beta|).$$
When $x_i$ is a good approximation of $x$, the value of $|\beta|$ is small. If we use shift as
 $$\sigma = x_i^T A x_i =\lambda \alpha^2 + (z^T Az) \beta^2. $$

We can obtain
$$2\mu_2/\mu_1 = 2\mu_2(\lambda - \sigma) = 2\mu_2 [(1-\alpha^2)\lambda + z^TAz \beta^2]. $$
Usually, $\lambda$ is not the largest eigenvalue and it is closer to $\sigma$ more than other eigenvalues of $A$.
Therefore, we can assume that $z^TAz$ and $2\mu_2$ are not large constants. From $\|x_i\|=1$, we can get
$1-\alpha^2=\beta^2$. With these results, we can draw the following result from (\ref{iequal3})
$$ \varepsilon < \frac{|\alpha|}{1+|\alpha\beta|}.$$

This shows that when $x_i$ is a good approximation of $x$,
the convergence of Algorithm \ref{vpower} does not require a very small $\varepsilon$.
From (\ref{speed}), the convergence rate  of Algorithm \ref{vpower} can be expressed as
\begin{equation}\label{speed2}
\frac{\mbox{tan}\angle(\tilde{y},x)}{\mbox{tan}\angle(x_i,x)} \leq \frac{2 \mu_2 |\beta + \tilde{\beta}|}{\mu_1|\alpha + \tilde{\alpha}|}
\leq \frac{2 \mu_2}{\mu_1} \frac{|\beta| + \varepsilon}{|\alpha| - \varepsilon} \frac{|\alpha|}{|\beta|}.
\end{equation}
If $\varepsilon=0$, the convergence rate is decided by $\frac{2 \mu_2}{\mu_1}$. When we set $\sigma = x_i^T A x_i$,
we have $\frac{2 \mu_2}{\mu_1}=O(\beta^2)$. This means that Algorithm \ref{vpower} is cubic convergence,
because in this choice of $\sigma$, one step of Algorithm \ref{vpower} is one step of Rayleigh quotient iteration.
When $\varepsilon \neq 0$, the convergence rate is slowed down. But, the convergence be damaged
only when  $\frac{|\beta| + \varepsilon}{|\alpha| - \varepsilon}> \frac{1}{\beta^2} $.
If  $\frac{|\beta| + \varepsilon}{|\alpha| - \varepsilon}$ is not very large, the convergence rate is also decided by $\frac{2 \mu_2}{\mu_1}$.

Along the standard lines of the inverse power method, the difference between ${\tilde{y}}$ and $y$ is almost parallel to the eigenvector. When the sequence ${x_i}$ begin to converge to the eigenvector, the inexact method can maintain the convergence trend very well, with a moderate accuracy of inner iteration.

As we described in Section 2, one has to compute the mostleft eigenvalue. Therefore, we can use an approximation eigenvalue as the target $\sigma$. For instance, we use the matlab command "eigs" to compute the approximate eigenpair. When the convergence rate slows down, we renew the target $\sigma$. One, we can use the generalized minimal residual method (GMRES) to compute the approximate solution of $(A-\sigma I )y = x_i$. If we are now to replace $x_i$ by the residue $r$ in $(A-\sigma I )y = x_i$, the IIPM becomes a two dimensional residual iterative method. For the residual iterative method, in Ref. \cite{Jia} it was shown that the inexact method can mimic the exact method with accuracy $\varepsilon_1=10^{-4}$.

For a non-Hermite matrix, we can compute the consequent approximate eigenvector $x_{i+1}$ from the subspace $span\{x_i, \tilde{y}\}$, to ensure that $x_{i+1}$ is a better approximation than $x_i$. From this point of view, one could as well call this method the modified IIPM.

\section{Numerical results}
\label{s5}

In this section, we first compute some matrices from the Matrix Market
by the exact inverse power method and modified inexact inverse power method to illustrate the validity of the theory analysis of the new method.
Here the exact method refers to the method that computes the exact solution at the third step of Algorithm 2.
Then we compute the practical problems which are derived from the two models of SEO to illustrate the
practicability of the new method. The two well-known models are the stochastic ABC model and the stochastic Kuramoto model. The phase space in both cases is a 3-torus, $X=T^3$, and the vector fields defining the noise, $e_a$'s, correspond to the additive Gaussian white noise,
\begin{eqnarray}
e_{1} \equiv e_x = (1,0,0)^T, e_{2} \equiv e_{y} = (0,1,0)^T, e_{3} \equiv e_{z} = (0,0,1)^T,
\end{eqnarray}
in the standard global coordinates on a 3-torus.

The flow vector fields of the two models are  given respectively as,
$$
F_{ABC} = (A\mbox{sin} z + C\mbox{cos} y)e_x + (B\mbox{sin} x + A\mbox{cos} z)e_y + (C\mbox{sin} y + B\mbox{cos} x)e_z.
$$
$$
\begin{array}{l}
F_{Kur}=\\
 (\omega_x - K/4(\mbox{2sin}x + \mbox{sin}(x+y) +  \mbox{sin}(x+y+z) - \mbox{sin}y - \mbox{sin}(y+z) ) )e_x\\
 + (\omega_y - K/4(\mbox{2sin}y + \mbox{sin}(x+y) + \mbox{sin}(y+z) - \mbox{sin}x - \mbox{sin}z))e_y \\
 +  (\omega_z - K/4(\mbox{2sin}z +  \mbox{sin}(y+z) + \mbox{sin}(x+y+z) - \mbox{sin}y - \mbox{sin}(x+y) ))e_z.
\end{array}
$$

A few remarks about the two models of interest are in order. First, the stochastic ABC model is a toy model for studies of astrophysical phenomenon of kinematic dynamo, \emph{i.e.}, the phenomenon of the generation of magnetic field by ionized flow of matter. As it was shown in Ref.\cite{Igor3}, the stochastic evolution of non-supersymmetric 2-forms of the STS of the stochastic ABC model,
\begin{eqnarray}
\partial_t \psi^{(2)} = -\hat H^{(2)}\psi^{(2)},
\end{eqnarray}
is equivalent to the dynamical equation of the magnetic field, $B$, in the kinematic dynamo theory,
\begin{equation}
\partial_t B = \partial \times F \times B + R_m^{-1} \bigtriangleup B,
\end{equation}
where $R_m=\Theta^{-1}$ is the inverse temperature known in the kinematic dynamo theory as the magnetic Reynolds number, and $\times$ denotes the standard vector product.

As to the Kuramoto model, it can be thought of as a model of coupled phase oscillators. This model also has many interesting applications. In particular, it may serve as a testbed for the studies of the phenomenon of synchronization that has attracted interest of scientists in biological \cite{Tass}, chemical \cite{Kiss}, physical\cite{Wiesenfeld} and other dynamical systems. An explicitly supersymmetric numerical representation of SEO on a square lattice of a 3-torus was proposed and described in the Appendix of Ref.\cite{Ours}.
This is the representation that we use in this paper.

All the experiments are run on  an Inspur Yitan NF5288 workstation with Intel(R)Core(TM)i5-3470s CPU 2.9GHz,
RAM 4G using Matlab R2012b  under the Linux system.

\begin{exm}
We compare the convergence between inexact method and exact method using a few examples. For comparison, we choose some matrices from the Matrix Market which the corresponding  Eq.(\ref{equ2}) can be directly solved.
(a) The first matrix is $H = A+A^T$, where $A$ is the matrix "rw5151".
(b) The second matrix is "cry10000".
(c) The third matrix is "bcsstk29".
\end{exm}
For each matrix, we use "exact" method and "inexact" method to compute the largest and smallest eigenvalues.
For the "exact" method, we require the solution to satisfy the $\epsilon_{\rm mach}=10^{-16}$. For the "inexact" method,
the accuracy of the inner iteration is $10^{-2}$.
The convergence process are shown in the following figures.
In the figures, "la" is the largest eigenvalue, "sa" is the smallest eigenvalue.

\begin{figure}[!h]
\begin{centering}
\includegraphics[height=4cm,width=4.8cm]{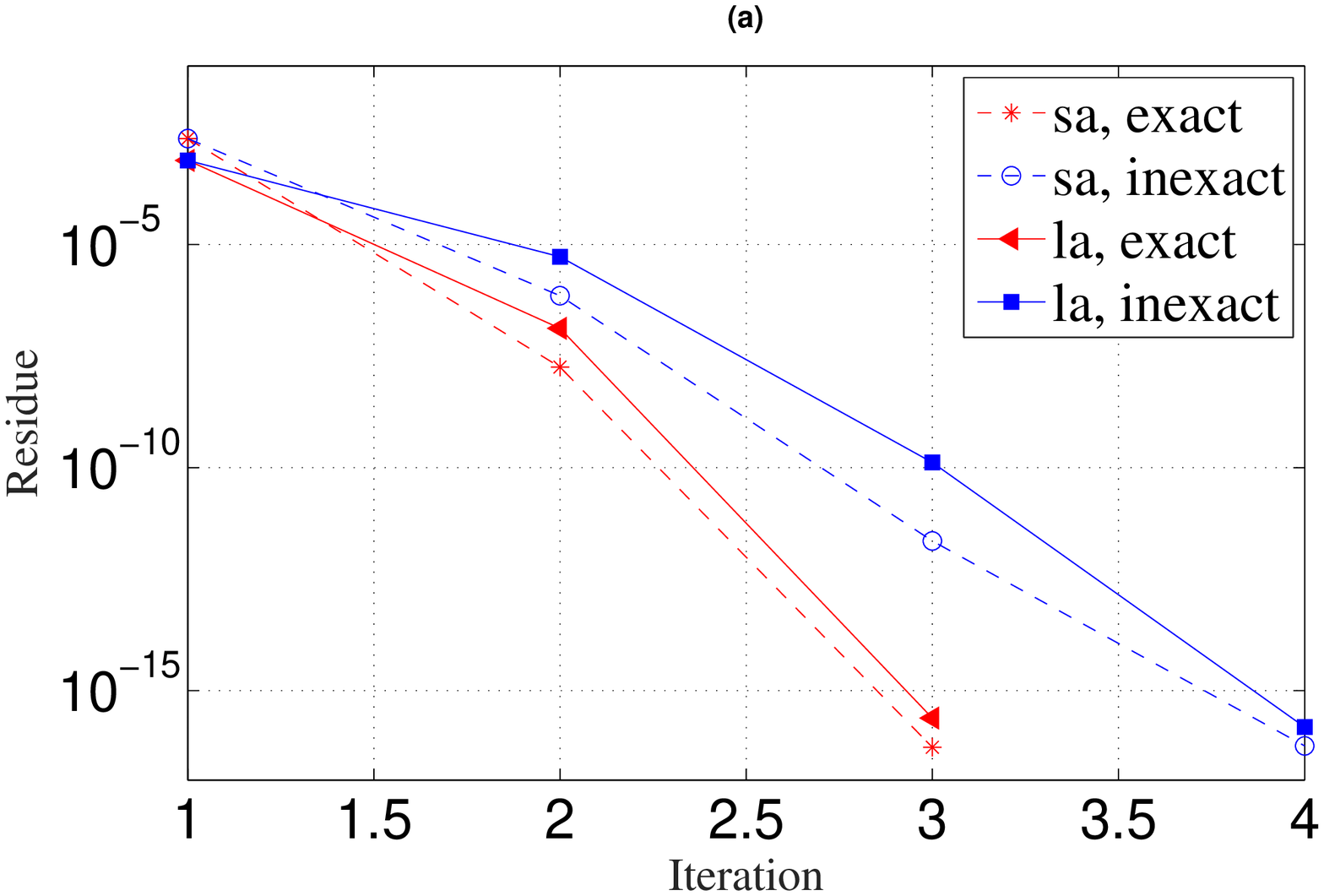}
\includegraphics[height=4cm,width=4.8cm]{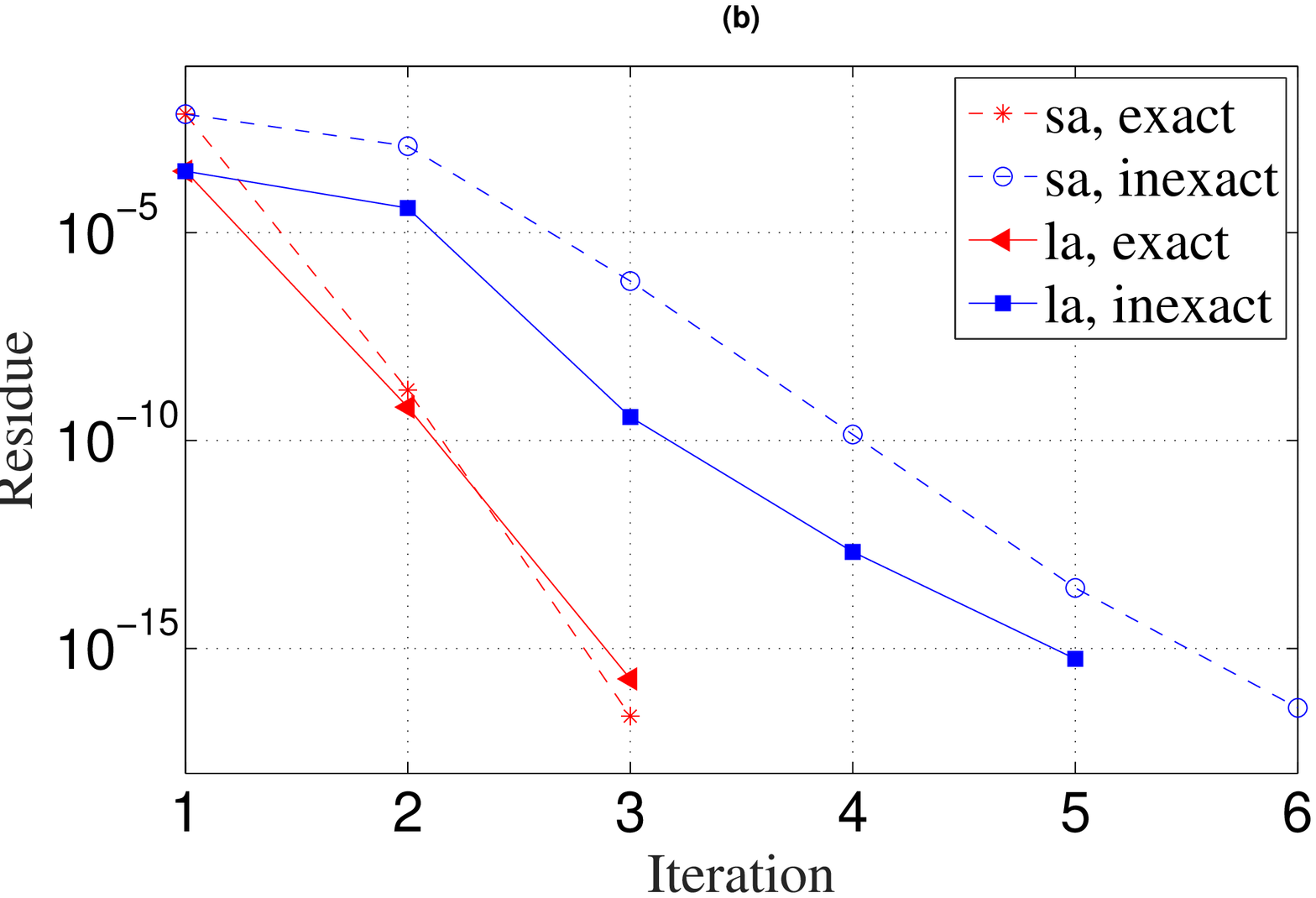}
\includegraphics[height=4cm,width=4.8cm]{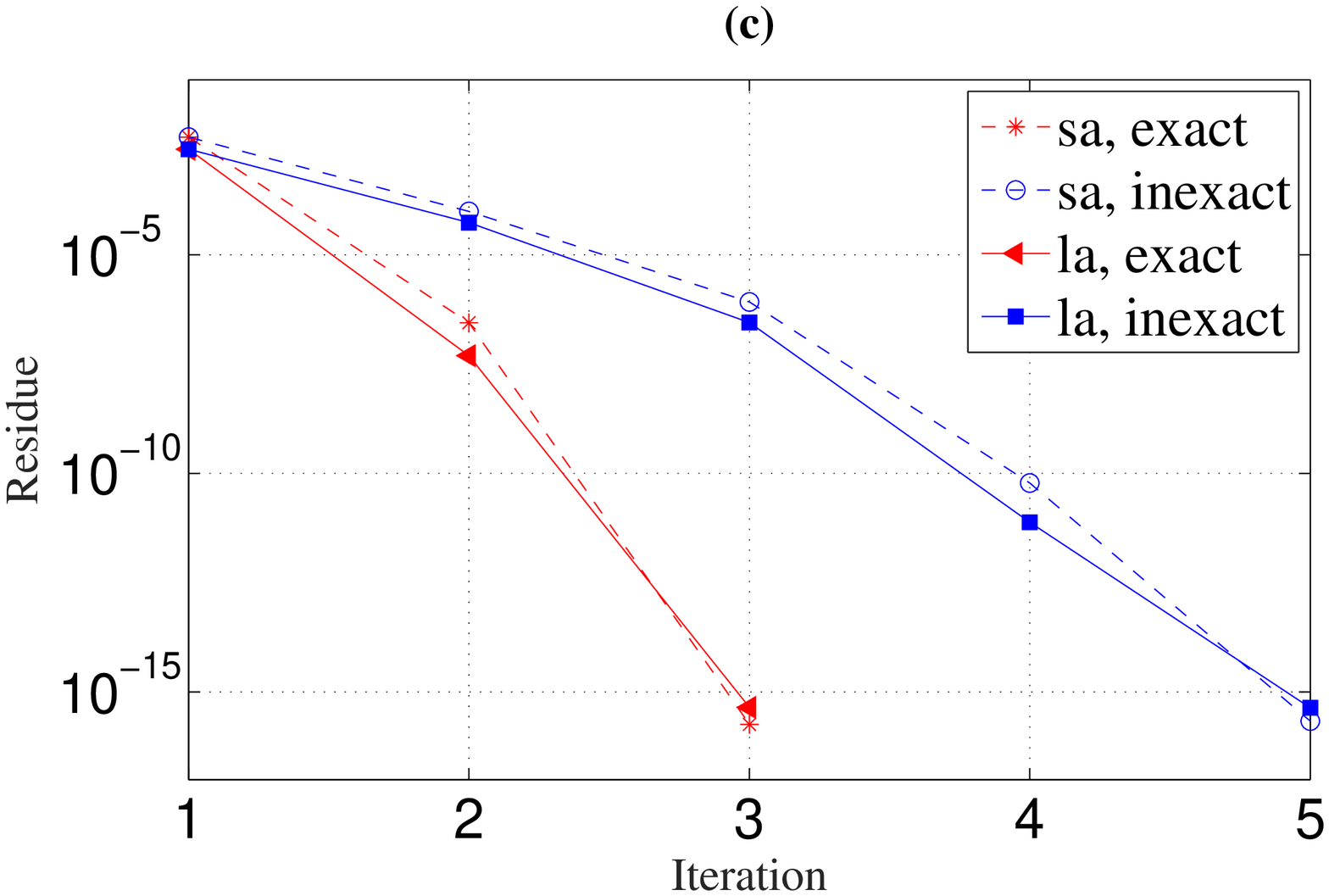}
\caption{The convergence of the exact (red curves) and inexact (blue curves) methods. The three subfigures represent convergence data, \emph{i.e.}, the norm of residue, of the largest (solid curves) and the smallest (dashed curves) eigenvalues as a function of the intereation number and for three different matrices (a-c). Even though the exact method converges faster than the inexact method in terms of the number of iteration, the inexact method demonstrates relatively good convergence. Furthermore, inasmuch each iteration of the inexact method is considerably faster than that of the exact method, the inexact method is actually much faster in terms of the real time of the computations.}
\label{figure1}
\end{centering}
\end{figure}
We can see from the figures that both methods converge quickly and smoothly.
The inexact method mimics the exact very well and it uses no more than three outer iterations.
The results confirm our theory and indicate that we can use it to solve more large problems.

In our application, there are several hundreds matrices need to be computed.
All matrices are too large to be solved using the exact method. Therefore, we use the new method to solve them.
Different matrix requires varies iteration number and cputime.
So, we give the information including their maximum, minimum, median of cputime and iteration number and so on.
\begin{exm}
In this example, we study the ABC model in the region $x, y, z\in [-\pi,\pi]$.
In order to study the influence of parameters $C$ and $R_m$, we select some points in the plane of $R_m$ and $C$.
For each pair of $(R_m, C)$, we discretize the ABC model to a matrix eigenvalue problem and analyze the system by the leftmost eigenvalue of the matrices.
\end{exm}
The points in the  $R_m$ and $C$ plane are $R_m = [1:1:14]$ and $C=[0.4:0.025:1.125]$.
For each point, the size of the matrix is 192000.
We compute the leftmost eigenvalue of the 420 large scale matrices. The real part of the eigenvalues are plotted in the following figure.
Where the circles represent the value less than or equal to zero and the plus represent the value greater than zero.
\begin{figure}[!h]
\begin{centering}
\includegraphics[height=6cm,width=8cm]{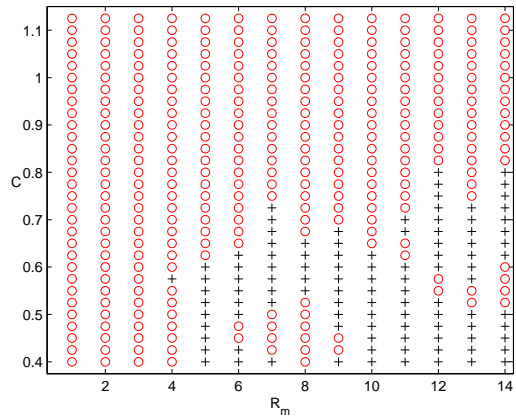}
\caption{Contour of ABC model. There are 420 points coresponding to 420 large scale matrices and the size of the matrix is 192000. We use the new method to compute the leftmost eigenpair of each matrix. The average cputime is 6.80s.
Where the circle and plus signs indicate the different state of eigenvalues. Based on this result, we can analyze the properties of the ABC model.}
\label{figure6}
\end{centering}
\end{figure}
In the computing process, we first use the matlab command "eigs" to compute an approximation of the leftmost eigenpair.
The convergence tolerance is $10^{-4}$. The Frobenius norm of the matrices are $O(10^3)$,
therefore, the absolute error of the approximate eigenvalues are $O(10^{-1})$. This is a modest request and "eigs" can
compute these result with a suitably large number of Lanczos vectors. But for most of the matrices, it is very hard to get more accurate results. In general, the desired eigenvalues are very close to the origin. So we first transform the leftmost eigenvalue to the module largest eigenvalue by a shift $\beta$.
Then we use "eigs" to compute the largest eigenvalue of matrix $A-\beta I$. Here the Frobenius norm of $A$ is a good choice of $\beta$.
We use the approximate eigenpairs as the target and starting vector of the inexact inverse  power method to compute the more accurate results.
The convergence tolerance of the outer interation is $10^{-10}$ and the maximum   iteration number of the outer iteration is 25.
We use GMRES to solve the inner linear systems with convergence tolerance $10^{-3}$.
We show the iteration number of outer and inner iteration in Table \ref{T-a},
\begin{table}[!h]
\begin{center}
\caption{Example 1, iteration number.}\label{T-a}
\begin{tabular}{|c|c|c|c|c|c|}\hline
{Iteration number} &Average & Maximum & Minimum &  Median & Total\\\hline
Outer & 1.07 &11 & 1 & 1 & 451 \\\hline
Inner & 1050.3 &4992 & 275 & 904 &441126 \\\hline
\end{tabular}
\end{center}
\end{table}
where "Outer" represents the iteration number of inverse power and "Inner" represents the number of Lanczos vectors in GMRES.
"Total" is the sum of all the 420 matrices. We also show the maximum, minimum, average and median of the iteration number.
We can see from the table that we used 451 times inverse power iteration and 441126 Lanczos vectors of GMRES to obtain all the desired eigenvalues.
The average outer and inner iteration number of all matrices are 1.07 and 1050.30 respectively.
The median  number of outer and inner iteration are 1 and 904.
The matrix corresponding to $R_m=10, C=1$  needs the most outer iterations. The matrices corresponding to $R_m=11, C=0.85$ and
$R_m=1, C=1.075$ need most and lest inner iterations respectively.

The cputime of all 420 matrices computation is 4576.34 seconds and the details of each part are shown in Table \ref{T-b}.
\begin{table}[!h]
\begin{center}
\caption{Example 1, details of cputime. }
\label{T-b}
\begin{tabular}{|c|c|c|c|c|c|c|}\hline
{cputime } & Average & Maximum & Minimum &  Median & Total  \\\hline
   eigs    & 6.80 &20.96 &2.51 & 6.67 & 2856.27  \\\hline
   GMRES   & 697.14 &13215.794 & 61.78 &529.22 & 292796.97  \\\hline
   Entire  & 705.39 &13223.85 & 66.19 & 537.62 & 296265.35  \\\hline
\end{tabular}
\end{center}
\end{table}

The time for computing the approximate eigenpairs is 2856.27 seconds.
For a single matrix, the maximum is  20.96s, the minimum is 2.51s, the average is 6.80s and median is 6.67s.
We use a total of 292796.97 seconds to solve all inner equations, and the average of each matrix is 697.14 seconds.
The total time of each  matrix to compute the eigenvalue is 705.39 seconds on average.

These results indicate that the most time-consuming part is the inner iteration. Thus reducing the number of outer iteration is very important to improve the computational efficiency.

\begin{exm}
In this example we study the Kuramoto model in the region $x, y, z\in [-\pi,\pi]$.
We 	analyze the influence of parameters $K$ and $D$. We set some points in the plane of ($K$, $D$).
For each pair of $(K, D)$, We use 30 lattice sites in each of the three directions to discretize the model to a matrix eigenvalue problem. The size of the corresponding matrix $A$ is 128625.
We can analyze the property of the system by computing the eigenvalues of $A$ .
\end{exm}
We set the values of $K$ and $D$  to $D = [0.0015:0.0027:0.042]$ and $K =[0.12:0.02:0.7]$.
From the real part of the  leftmost eigenvalues of the 480
large scale matrices, we get the following figure,
where the circles represent the value less than or equal to zero and the plus represent the value greater than zero.
\begin{figure}[!h]
\begin{centering}
\includegraphics[height=6cm,width=8cm]{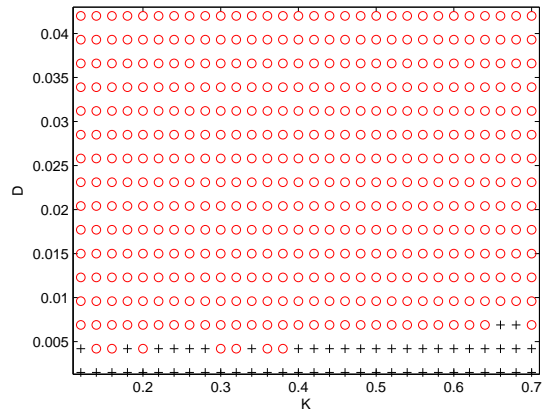}
\caption{Contour of Kuramoto model.
There are 480 points coresponding to 480 large scale matrices and the size of the matrix is 128625. We use the new method to compute the leftmost eigenpair of each matrix. The average cputime is 130.6s.}
\label{figure2}
\end{centering}
\end{figure}
We first use "eigs" to compute the approximation eigenpairs $(\tilde{\lambda},\tilde{x})$ satisfying
$\frac{||A\tilde{x} -\tilde{\lambda}\tilde{x}||}{||A||} \leq10^{-5}$. The order of $||A||$ for all matrices are $O(10^5)$, so the approximate eigenvalues hardly have any algebraic precision for the absolute error. Even so, this is a difficult task for eigs.
We use the approximate eigenpairs as the target and starting vector of the inexact inverse  power method to compute more accuracy results.
The convergence tolerance of the outer iteration is $10^{-11}$ and the maximum iteration number of the outer iteration is 25.
We use GMRES to solve the inner linear systems with convergence tolerance $10^{-4}$.
We show the iteration number of outer and inner iterations in the following table.
\begin{table}[h]
\begin{center}
\caption{Example 2, iteration number. }
\label{1-a}
\begin{tabular}{|c|c|c|c|c|c|}\hline
{Iteration number} &Average & Maximum & Minimum &  Median & Total\\\hline
Outer & 7.75 &47 & 1 & 5 & 3953\\\hline
Inner & 1283 &17523 & 5 & 288 & 654330\\\hline
\end{tabular}
\end{center}
\end{table}
Where "Outer" is the outer iteration number and "Inner" is the number of Lanczos vectors of GMRES.

The cputime of all 480 matrices is 220942.50 seconds and the details of each part are shown in the following table.
\begin{table}[!h]
\begin{center}
\caption{Example 2, details of cputime. }
\label{1-b}
\begin{tabular}{|c|c|c|c|c|c|c|}\hline
{cputime } & Average & Maximum & Minimum &  Median & Total  \\\hline
   eigs    & 130.60 &2175.40 & 7.84& 42.88 &  62687.04\\\hline
GMRES      & 323.67 &13343.00 & 0.49 & 33.03 & 155360.50  \\\hline
Entire     & 460.30 &13377 & 37.96 &  202.59& 220942.50  \\\hline
\end{tabular}
\end{center}
\end{table}
The time for computing the approximate eigenpairs is 62687.04 seconds.
For a single matrix, the maximum is 2175.40 seconds, the minimum is 7.84 seconds, the average is 130.60 seconds and median is 42.88 seconds.
We use a total of 155360.50 seconds to solve the inner equations, and the average of each matrix is 323.67 seconds.
The total time of each matrix to compute the eigenvalue is 460.30 seconds on average.

The results show that this model is more difficult than the previous model. It takes more external iterations and the inner iteration is still the most time-consuming part. The large difference between median and average of the inner cputime indicate that there is a big difference of the inner cputime between different matrices.

\section{Conclusion}
\label{conclusion}

In this paper, we proposed what we call the inexact inverse power method (IIPM) for numerical diagonalization of sparse matrices. This method allows to notably save computational resources as compared to its parental well-established inverse power method. We applied IIPM to the problem of the finding the ground state of the stochastic evolution operators of the stochastic ABC and Kuramoto models and our results demonstrate that IIPM provides solution at acceptable computational time in situations when IPM would fail if using only the resources of a typical desktop computer.






\bibliographystyle{model1-num-names}
\bibliography{sample.bib}



\end{document}